\documentclass[10pt,twoside,reqno]{amsart}

\usepackage{amsfonts}
\usepackage{amssymb}
\usepackage{amsmath}    %multi-line etc. Subsumes ams-text,amsbsy,amsopn
\usepackage{amsthm}     %provide proof env
\usepackage{amscd}      %commutative diagrams
\usepackage{euscript}
\usepackage[english]{babel}

{
       \newtheorem{theorem}{Theorem}[section]
       \newtheorem{proposition}[theorem]{Proposition}
       \newtheorem{lemma}[theorem]{Lemma}

{\theoremstyle{definition}
       \newtheorem{remark}[theorem]{Remark}
}

\newcommand{\RR}{{\mathbb{R}}}

\newcommand{\QQ}{{\mathbb{Q}}}

\newcommand{\Ric}{\operatorname{Ric}}
\newcommand{\R}{\operatorname{R}}
\newcommand{\Vol}{\operatorname{Vol}}
\newcommand{\tr}{\operatorname{tr}}

\newcommand{\Pic}{\operatorname{Pic}}
\newcommand{\Id}{\operatorname{I}}
\newcommand{\HH}{\operatorname{H}}
\newcommand{\mult}{\operatorname{mult}}
\newcommand{\p}{\partial}

\newcommand{\ddc}{\operatorname{d d^c}}

\newcommand{\bj}{\bar j}

\begin{document}

\title
[K\"ahler-Ricci flow on projective varieties]
{K\"ahler-Ricci flow and the Minimal Model Program for
projective varieties}
\author[P. Cascini and G. La Nave] {Paolo Cascini and Gabriele La Nave }
\address{Department of Mathematics, University of California at Santa
    Barbara, Santa Barbara, CA 93106, US}
\email{cascini@cims.nyu.edu}
\address{Lehigh University and Courant Institute, 14 E. Packer Street,
Bethlehem PA}
\email{gal204@lehigh.edu}

\begin{abstract}
In this note we propose to show that the K\"ahler-Ricci flow fits
naturally within the context of the Minimal Model Program for
projective varieties. In particular we show that the flow detects, in
finite time, the
contraction theorem of any extremal ray and we analyze the
singularities of the metric in the case of divisorial
contractions for varieties of general type. In case one has a smooth
minimal model of general type (i.e., the canonical bundle is nef and
big), we show infinite time existence and
analyze the singularities. 

\end{abstract}

\maketitle

\section {Introduction}

One of the most important problems in Algebraic Geometry is the
quest for a Minimal Model (i.e., a variety which is birationally
equivalent to the given one whose canonical bundle is {\it nef}),
dubbed ``the Minimal Model Program''. This entails the application
of a complicated algorithm, which has been proved to work in
dimension $3$ by the collaborative effort of Mori, Kawamata,
Koll{\'a}r, Shokurov et al. (see \cite{kmm} for a survey) and very recently in dimension $4$ \cite{HM}. In complex
dimension $2$ the theory is much simpler and is due to the Italian
school of algebraic geometry (like Enriques, Castelnuovo and
Severi) and put in a modern and more precise  framework by Zariski
and Kodaira. In case the variety is of {\it general type}  (i.e.,
varieties for which the space of holomorphic sections of $K_M ^k$
grows like $k^n$, where $n=\dim M$), then out of a minimal model
one can produce the so called {\it canonical model}, i.e., a
birational model whose canonical bundle is ample, or, in other words, with negative first
Chern class $c_1$. 

On the other hand, around the 1980's, building on the foundational
work of Hamilton in the Riemannian case, H. D. Cao studied the
K\"ahler-Ricci flow for canonical metrics on manifolds with definite
first Chern class.
$c_1$, reproving (see \cite{ca}) in particular Calabi's conjecture and the existence
of K\"ahler-Einstein metrics in case $c_1 <0$ (the solution of the
conjecture is originally due to S.T. Yau, cf. \cite{y1}). Cast in an
algebro-geometric light, this last result says that any smooth
    projective variety with ample canonical bundle admits a
    K\"ahler-Einstein metric (see \cite{DK} if $M$ is a projective
    variety with orbifold singularities).

In this note we propose to draw a connection between the two theories for projective
varieties of general type, and in fact prove that in complex dimension two the K\"ahler-Ricci flow
produces the {\it canonical model},
generalizing Cao's result.

In order to state the two main theorems we prove, we need to fix some
notation.
Let $M$ be a projective variety and $K_M$ its {\it canonical} line
bundle. If $K_M$ is not {\it nef}, there exists a (complex) curve   $C$  in $M$ on which $K_M \cdot C
    <0$. By the
rationality theorem
 (e.g. see \cite{kmm}), there exists a  nef line
bundle $L$ such that $A= L-r K_M$ is ample, for some rational number $r>0$. Therefore (by the
base-point-free theorem) $L=A+rK_M$ is
{\it semiample}, that is to say some power $L^n$ of $L$
defines a holomorphic map $c:M \to M'$. Moreover the morphism $c$
contracts only the curves that are homologically equivalent to $C$.

In this context we prove:

\begin{theorem}\label{main1}
Let $M,M'$ be as above,  and let $g_0$ be a metric
in the ample class $A$.
Then the K\"ahler-Ricci flow:
$$\frac {\p g_{i \bj}}{\p t} = -\Ric _{i\bj}-g_{i \bj}$$
flows the metric in the class $c_1\left(A+ a(t)(K_M -A )\right)$ for
$a(t)=1-e^{-t}$ and it develops a singularity at (the finite
time) $T=\log (r+1)$. 

Furthermore, if $M$ is of general type, then the
    singular locus $S$ of $g(T)$ is contained in a proper subvariety
    of $M$. If, in addition, $M'$ is smooth (in particular $c$ is a
    divisorial contraction), then $g(T)$ induces a smooth metric on
    $M'$. 
\end{theorem}

In case $M$ is a smooth minimal model (i.e. $K_M$ is nef), we prove:

\begin{theorem}\label{main2}

Let $M$ be a smooth projective manifold with $K_M$ big and nef,
then the normalized K\"ahler-Ricci flow as above exists for all
time. Moreover, if the canonical model $M'$ of $M$ admits only
    orbifold singularities, then the limit for $t\to + \infty$ of
    $g(t)$ is locally equal to $g_0 +d d^c v$ where $v$ is a bounded function which is smooth away
from an analytic subvariety.
\end{theorem}

In particular if $M$ is a projective surface of general type, then all
    the hypotesis of the previous two theorems are satisfied, and
    therefore, after a finite number of steps, the K\"ahler-Ricci flow
    converges to a K\"ahler-Einstein metric on  the canonical model of $M$. 
 
A more general case will be studied in a forthcoming paper. 

\vskip .5 cm 
Let us remark that Theorem \ref{main2} was a conjecture of Tian and
that it has been proved independently (with weaker assumptions) by G. Tian and
Z. Zhang (cf. \cite{tianzhang}).

The techniques involved are mainly the reduction of the K\"ahler-Ricci
flow at hand to a scalar parabolic PDE
similar to Cao's (\cite{ca}), and the coarse study of its
singularities, in the context of orbifolds.

We are grateful to F. Bogomolov, H.D. Cao, J. M\textsuperscript cKernan and M. McQuillan for many
    valuable discussions.

\vskip .5 cm

{\em Notations and Conventions.}
 On a $n$-dimensional complex manifold
    $M$, with local holomorphic coordinates $z_1,\cdots, z_n$, we will write
    $d=\p +\bar \p$ and $d^c= \frac {\sqrt
    {-1}}{4\pi} (\p - \bar \p)$, so that $\ddc = \frac {\sqrt
    {-1}}{2\pi} \p\bar\p$.
     Given a K\"ahler metric   $g=\frac
    {\sqrt{-1}}{2\pi} \sum g_{i\bj}dz_i\wedge dz_{\bj}$, we will
    denote by $\Ric_{g}$ its Ricci curvature, locally given by
    $\Ric_g=-\ddc \log \det (g_{ij})$.
    We will often use matricial notation for $(1,1)$ form, e.g. $g^{-1}$
     locally denotes the inverse matrix $g^{i\bar j}$ of $g$. As
     usual, given two quadratic forms on a given vector space, $q_1$ and $q_2$, with $q_1$
     non-singular, we set $tr_{q_1} q_2= tr (q_1 ^{-1}q_2)$. 
\vskip 1 cm
\section{The scalar equation}\label{scalar.equation}

Let $(M,g_0)$ be a projective manifold with a K\"ahler metric
    $g_0$ belonging to the ample class  $A$.
The first thing to observe is that if the K\"ahler metric $g(t)$ satisfies
    the K\"ahler-Ricci flow,
\begin{equation}\label{ricciflow}
\left\{
\begin{aligned}
& g'(t)= -\Ric_{g(t)} - g(t)\\
& g(0)=g_0
\end{aligned}
\right .
\end{equation}
 then, $g(t)$ is in the class of $c_1 (A(t))$ where $A(0)=A$ and
    $A(t)$ is an ample class on $M$.

By abuse of notation, and since $-\Ric_g\in
    c_1(K_M)$ for any K\"ahler metric $g$, we have

$$\left\{
\begin{aligned}
&\p_t A(t) = K_M - A(t)\\
&A(0)=A
\end{aligned}\right.$$
or, therefore, by solving the ODE, we have the following
\begin{proposition}\label{class.eq}
If $g(t)$ is a solution for the K\"ahler-Ricci flow (\ref{ricciflow}),
    then $g(t)$ belongs to the class
$$A(t)=A+a(t)(K_M-A)$$
with $a(t)=1-e^{-t}$.
\end{proposition}

With this at hand, let us proceed to reduce our equation to a scalar
one.

Let $\eta_0=-g_0-\Ric_{g_0}$.
In particular, $\eta_0$ is an element in the  class $K_M -A$.
Thus, given any element $\eta$ in the same class (we will allow ourselves the
    freedom of choosing $\eta$ later on, depending on the situation), there exists a
    smooth function $f=f_\eta$ on $M$, such that:
$$\eta_0=\eta+\ddc f.$$

Moreover, if we consider the $(1,1)$-form 
\begin{equation} \label{g0}
g_0(t)=g_0+a(t) \eta,
\end{equation} 
then $g_0(t)$ belongs to the class $A(t)$, and therefore 
    there exists a smooth function $u$ (also depending on $\eta$),
    such that

$$g(t)=g_0(t) + \ddc u.$$

Thus differentiating, and using the fact that $\Ric_{g(t)}=-\ddc \log
    \big(\det g(t)\big)$, equation (\ref{ricciflow}) becomes 
$$a'(t) \; \eta + \ddc(\p_{t}u)=
\ddc \log ~ \det\big( \frac{ g_0(t)  +\ddc u}{g_0}\big)
+ \eta_0 - a(t)\eta -\ddc u$$
and therefore:

$$\ddc\left(\p_{t}u)\right)= \ddc \log~ \det\left(
    \frac{g_0+a(t)\eta+\ddc  u}{g_0}\right)- \ddc u +\ddc f.$$
Thus, by the $\p \overline {\p}$-lemma, there exists a smooth function $\phi(t)$
(depending only on $t$)  such that:

$$\p_t u= \log~ \det\left(
    \frac{g_0+a(t)\eta+\ddc u}{g_0}\right)-u +f +\phi.$$

Moreover $\phi$ satisfies:
$$\int _M \exp\left(\p_{t}u+u-f\right)~dV_0 =\Vol(M) e^{\phi(t)}$$
where $dV_0$ is the volume form of the metric $g_0$, and $\Vol(M)$ is
    the volume of $M$ with respect to the metric $g(t)$.

Renormalizing $u$ so that $\phi=0$, we have the following:

\begin{lemma}
The K\"ahler-Ricci flow as in theorem \ref{main1} is equivalent to the
scalar equation:
\begin{equation}\label{scalar}
\left\{ \begin{aligned} &\p_{t} u= \log ~ \det\left(
    \frac{g_0+a(t)\eta+\ddc u}{g_0} \right)-u +f\\
&  u(x,0)=0
\end{aligned}\right.
\end{equation}
\end{lemma}

Notice that the solution $u$ to eq. (\ref{scalar}) is dependent on the
choice of $\eta$, but $g$ is not. In fact we have:

\begin{lemma}\label{oneone}
There exists a one-to-one correspondence between the solutions $u$ of
(\ref{scalar}) with pair $(\eta,f)$ and the solutions $u'$ with pair $(\eta ',
f')=(\eta ' , f')=(\eta -dd^c h, f+h)$, given by $u'=u+a(t)h$. In
    particular $h$ does
    not depend on $t$.
\end{lemma}
\vskip 1cm

\section {Maximal existence time Case I:  $K_M$ not nef}\label{longtime}

Given a K\"ahler manifold $(M,g_0)$, we first investigate the behavior of the
solutions to equation (\ref{scalar}), in the case $M$ is a projective
    variety with canonical bundle not nef.

The equation is parabolic and general theory implies short time
existence (e.g. see \cite{ta}). Moreover we have:

\begin{proposition}\label{c0bound}
Suppose $u$ is a solution to eq. (\ref{scalar}) in $M\times
(0,t_0)$, for some time $t_0>0$. Then there
exists a uniform constant $C>0$ such that:
$$ |u| <C.$$
In particular, if the flow exists for any $t<t_0$, then for any
sequence $\{t_i\} \subset [0,t_0)$, we have that 
$\lim _{t_i\to T}u(x,t_i)$ is continuous (up to taking a sub-sequence).
\end{proposition}

\vskip .5 cm

In order to prove the proposition, we have to choose a suitable
    $\eta$. Although $u$ depends on $\eta$, from lemma \ref{oneone} it follows that
    the fact that $u$ is
    bounded does not depend on that particular choice.

\vskip .5 cm

As stated  in the introduction, by the rationality theorem and the
    base point free theorem, since $K_M$ is not nef, there exists a rational number $r$ such
    that $L=A+rK_M$ is semi-ample and defines a contraction $c:M\to
    M'$ of an extremal ray.
In fact,
\begin{equation}\label{def.r}
r=\max \{s\in \mathbb Q|A+sK_M \text{ is nef}\}.
\end{equation}

In particular, by proposition \ref{class.eq}, it follows that  there
    can be a solution for (\ref{ricciflow}) at most up to time
    $$T=\log(r+1).$$ 
In fact, the class $A(T)=(r+1)^{-1}L$ is not ample
    and therefore it cannot contain any metric.

On the other hand, $L$ is semi-ample,  i.e. it is the pull-back of some
    ample class $A'$ on $M'$ with respect to the contraction map
    $c$. Therefore there exists a  non-negative $(1,1)$ form
    $\eta_L$ in $c_1(L)$.

By prop. \ref{class.eq}, we can write $$A(t)=\frac 1 r (a(t)L+b(t)A)$$
where $a(t)=1-e^{-t}$ and $b(t)=(r+1)e^{-t}-1$. in particular we can choose $\eta$ as
\begin{equation}\label{eta.etal}
\eta= \frac 1 r(\eta_L-(r+1)g_0),
\end{equation}
so that the $(1,1)$ form $g_0(t)$, defined in (\ref{g0}), is given by 
 $$g_0(t)=g_0+a(t)\eta=\frac 1 r (a(t)\eta_L + b(t) g_0),$$
and therefore it is a metric for any $t<T$. 

\vskip .5 cm

We need the following

\begin{lemma}\label{supersolution}
There exists a bounded super-solution $u^+$ (resp. a bounded
    sub-solution $u^-$) for (\ref{scalar}),
    depending only on the time $t$, defined for any $t\in [0,T)$,
    i.e.
$$\begin{aligned}
& \partial_t {u^+} + u^+ \ge
\log(\det(g_0(t) /{g_0})) + f \\
(~\text{ resp. }\quad &\partial_t {u^-} + u^- \le
\log(\det(g_0(t) /{g_0})) + f\quad).
\end{aligned}
$$
\end{lemma}
\vskip .5 cm

\begin{proof}
Since $M$ is compact, and since $f$ is a smooth function on $M$ and $a(t)$ and $b(t)$ are bounded
    functions on $[0,T)$,  there exists a positive constant $K$ such that
$\log(\det(g_0(t) /{g_0})) + f<K$.

Therefore, in order to define a super-solution for (\ref{scalar}), it is enough
    to choose $u^+$ as a solution of
$\partial_t {u^+} + u^+=K$, with $u^+(0)=0$, i.e.
$$u^+ = (1-e^{-t})K.$$

On the other hand,  $a(t), b(t)$ are non-negative  in the interval
    $[0,T)$ and   since $\eta_L$ is semi-positive, it follows
$$\det(g_0(t)/g_0)=\det\left( \frac 1 r \cdot \frac {a(t)\eta_L + b(t)
    g_0}{g_0}\right)\ge \left(\frac {b(t)} r \right )^n.$$

It can be easily checked that $\displaystyle \int_0^T\log~b(t)~dt>-\infty$
    and therefore, given $K$ such that $f>K$, we can define $u^-$ as
    the solution for
$$\left \{
\begin{aligned}&\p_{t} u^- + u^- = \log \left(b(t)/ r \right )^n + K
    \\
&u^-(0)=0
\end{aligned}
\right .
$$
Thus $u^-$ is a bounded sub-solution
    for (\ref{scalar}).
\end{proof}

\vskip 1 cm
\begin{proof}[Proof of Proposition \ref{c0bound}]
\

Let $u^-$ and $u^+$ be as in lemma \ref{supersolution}.
By the comparison principle, we will show that for any solution $u(t)$
    of (\ref{scalar}), we  have
$$u^-(t)\le u(t)\le u^+(t).$$
Since $u^-$ and $u^+$ are bounded function, the proposition will follow.

Let $w=u-u^-$. Then, since $u^+$ is a super-solution for
    (\ref{scalar}) that depends only on the time $t$, we have that for
    any $t<T$,
$$\partial_t {w} + w \le
\log(\det(g_0(t)+\ddc w) /{g_0(t)}).$$

Therefore if $\displaystyle \overline w(t)=\max_M (w(\cdot, t))$, then
$\partial_t {\overline w} + \overline w \le 0$. Since $w(\cdot, 0)=0$, it
    follows $w\le 0$, i.e. $u\le u^+$. Similarly, it follows $u\ge
    u^-$.
\end{proof}

\vskip 1 cm

\begin{remark}Let us define $$F(x,t,r,p,X)=e^{p+r-f(x)}-\det((g_0(t)+X)/g_0).$$
Then the scalar equation (\ref{scalar}) for the K\"ahler-Ricci flow  is equivalent to the equation
\begin{equation}\label{viscosity}
\left\{
\begin{aligned}
&F(x,t,u(x,t),\p_t u(x,t), \ddc u(x,t))=0\\
&u(x,0)=0
\end{aligned}
\right .
\end{equation}
It is easy to check that the operator $F$ is {\it proper},
(according to the definition (0.2) and (0.3) in \cite{CIL}). 

Therefore by Perron's method (e.g., see theorem 4.1 in \cite{CIL}), 
the existence of a super-solution and a sub-solution for
    (\ref{viscosity}), guaranteed by
    lemma \ref{supersolution}, implies the existence of  a
weak solution (in the sense of viscosity solutions) for (\ref{viscosity}), for any $t<T$.

On the other hand, we are going to show that equation (\ref{scalar})
    admits a strong solution $u$ at any time $t<T$.  In particular,
    from the boundness of $\p_t u$ it will follow that  $g(t)$ is a K\"ahler metric
    for any $t<T$.
\end{remark}
\vskip .5 cm

We will need the following
\begin{lemma}\label{trace}
Let $\psi$ be a non-negative $(1,1)$-form in an ample
    class  on $M$. Given any $(1,1)$-form $\eta$, there exists $\eta'$
    (resp. $\eta''$) in the same class as $\eta$ and a constant $C$ such that
$$\tr_\psi\eta'>C \quad (\text{resp. }\tr_\psi\eta''<C ) $$
where it is defined.
\end{lemma}
\begin{proof}
Since $\psi$ belongs to an ample class, there
    exists a positive $(1,1)$-form $\psi_0$ cohomologous to $\psi$,
    i.e. $\psi_0=\psi + \ddc \phi$ for some $\phi\in
    C^\infty(M)$.

Thus, we have
$$\begin{aligned}
\tr_\psi\eta&=\tr ~(\psi^{-1} - \psi_0^{-1})\cdot \eta + \tr_{\psi_0}\eta \\
&=\tr ~(\psi^{-1}\cdot\ddc \phi\cdot \psi_0^{-1}\cdot\eta) + \tr_{\psi_0}\eta.
\end{aligned}
$$

Given a constant $a$, we can define $\eta'=\eta + a \ddc\phi$,
    that is a $(1,1)-$form cohomologous to $\eta$, such that
$$\begin{aligned}
\tr_\psi\eta'&=\tr_\psi \eta + a\tr_{\psi}\ddc \phi \\
&=\tr ~(\psi^{-1}\cdot\ddc\phi)\cdot (\psi_0^{-1}\cdot\eta + a \Id) + \tr_{\psi_0}\eta.
\end{aligned}$$
By the positivity of $\psi_0$, it follows that in the locus and in the
    directions where $\psi$ is zero, we have that $\ddc \phi$ is
    positive. Therefore, if $a$ is large enough so that
    $\psi_0^{-1}\eta+a \Id>0$, by the compactness of $M$, there exists a
    constant $C$ such that $\tr_\psi\eta'>C$.

Similarly, by choosing $\eta''=\eta-b\ddc\phi$ with $b$ large
    enough, there exists a constant  $C$ such that $\tr_\psi\eta''<C$.
\end{proof}

\vskip .5 cm

\begin{lemma}\label{rescale}
By rescaling the initial metric $g_0$, by a positive constant $K$, the
    singular locus for the solution $g(t)$ of the K\"ahler-Ricci flow
    (\ref{ricciflow}) at maximal time $T$, does not change.
\end{lemma}

\begin{proof}
Let $g(t)$ be a solution for the K\"ahler-Ricci flow with initial
    metric $g_0$, and let $K>0$.

If $\tilde g(s)=k(s) g(t(s))$, with
$$k(s)=(K-1)e^{-s}+1 \qquad \text{and} \qquad t(s)=\log\left(\frac
    {e^s+K-1} K\right ),$$
then $\tilde g(s)$ is a solution for the rescaled K\"ahler-Ricci flow
\begin{equation}
\left\{
\begin{aligned}
& \tilde g'(s)= -\Ric_{\tilde g(s)} - \tilde g(s)\\
& \tilde g(0)=Kg_0
\end{aligned}
\right .
\end{equation}
In particular, we have that the singular locus of $g(\log(r+1))$ coincides with the
    singular locus of $\tilde g(\log(Kr+1))$.
\end{proof}

\vskip .5 cm

We can now prove:

\begin{proposition}\label{ut.bound}
For any $t_0\in(0,T)$  there exist constants $C_0$, $C$, with
    $C$ independent of $t_0$, such that, for
    any $t<t_0$ and as long as there exists a solution for
    (\ref{scalar}),  we have
$$C_0<\p_{t} u <C $$
\end{proposition}

\begin{proof}
For ease of notation, let us denote $v=\p_t u$.
By taking the derivative
of
    (\ref{scalar}) with respect to $t$, we have that $v$ is a solution
    of
\begin{equation}\label{heat.scalar}
\left\{
\begin{aligned}
\p_t v &= \Delta_{g(t)} v + a'(t) \tr_{g(t)} \eta - v\\
v(0) &= f
\end{aligned}\right.
\end{equation}
where $\Delta_{g(t)}$ denotes the Laplacian with respect to the metric
    $g(t)$.

Let us first show that there is a uniform  upper bound for $v$. By rescaling $A$
if necessary, by lemma \ref{rescale} we can suppose
    without loss of generality that $A-K_M$ is ample. Thus we can
    choose $\eta$ as a negative $(1,1)$-form in the class of $K_M-A$
    and, since $v$ is a solution of (\ref{heat.scalar}), it follows
    that for this choice of $\eta$, we have
$$\p_t v \ge \Delta_{g(t)} v -v .$$
Thus by the Maximum Principle, $v$ is uniformly bounded from above.

Let us suppose now that there exists $t_0\in (0,T)$, such that, for any
    $t<t_0$, there exists a solution for (\ref{scalar}) , but
    $\inf_M v(\cdot, t_0)=-\infty$.

Since $t_0<T$, the non-negative $(1,1)-$form $\psi=g(t_0)$ belongs to the
    ample class $A(t_0)$. Therefore, by lemma
    \ref{trace}, there exists a $(1,1)-$form $\eta'$ in the class of
    $K_M-A$, and a constant $C$ such that $\tr_{\psi}
    \eta'>C$. Moreover, since $g(t)$ is positive for any $t<t_0$, we
    can choose $C$ so that $\tr_{g(t)}\eta'>C$ for any $t\le t_0$.

We can suppose, without loss of generality, that $\eta=\eta'$. In
    fact, by lemma \ref{oneone} the solution $v'$ for (\ref{heat.scalar}), associated to
    $\eta'$ will differ from $v$ by $v'=v-a'(t)h$, for some
    $C^\infty$ function $h$ on $M$, not depending on $t$. Thus, the locus where the
    limit of $\inf_M v$ and of $\inf _M v'$ are $-\infty$ as $t\to t_0$, coincide.

Thus we have
$$\partial_t v \ge \Delta_{g(t)} v - v + a'(t)C
$$
and by the Minimum Principle, it follows
$$\inf_M v(\cdot, t_0)\ge e^{-t_0}(Ct_0+\min_M f)>-\infty,$$
which gives a contradiction, so $v$ must be bounded from below on
$M\times [0,t]$ for any $t<T$ .
\end{proof}

\vskip 1 cm

We can now proceed to the $C^1$ and $C^2$ estimates. Let us denote by
    $\Delta=\Delta_{g_0(t)}$, the Laplacian with respect to the
    metric $g_0(t)$.  Moreover let us denote by $\R_{i \bar i j \bar j}(g_0(t))$
    the bisectional curvature of $g_0(t)$ (as mentioned above, we can
    choose $\eta$ as in (\ref{eta.etal}) so that $g_0(t)= g_0 +a(t)
    \eta$ is a metric for any $t<T$).

\begin{lemma} \label{laplacianestimate}
Let $t_0\in (0,T)$. There exists a constant $C_0$,  such that for any $t<t_0$ and as long
    as there exists a solution for (\ref{scalar}), we have
$$|\Delta u| \leq  C_0.$$
\end{lemma}
\begin{proof}
Since $g(t)=g_0(t) + \ddc u$,
as long as there exists a solution $g(t)$, we
    must have 
$$0<tr_{g_0(t)}g(t)=n + \Delta u$$ 
and
    therefore the lower bound follows immediately from the fact $g(t)$
    is a metric.

Let
\begin{equation}\label{def.F}
F:= \p_t u + u -f - \log(\det(g_0(t)/g_0)).
\end{equation}
   Then
$$ \det \left( g_0(t) + \ddc u\right)= e^F\det (g_0(t)).$$
From (2.22) in \cite{y1}, the following inequality holds for any
    $\lambda\in \RR$:
$$\begin{aligned}
 \Delta_{g(t)} \left(e^{-\lambda u}  (n+\Delta u)  \right)
    \geq  e^{-\lambda u}   \big( \Delta F& - n^2\inf \R_{i\bar
    {i} j \bj}(g_0(t))\\  &+  e^{-\frac F {n-1}}(\inf \R_{i\bar {i} j
    \bj}(g_0(t)) + \lambda)
(n+ \Delta u)^{\frac{n}{n-1}} -\lambda n(n+\Delta u)\big).
\end{aligned}
$$
From now on, we choose $\lambda$ to be any
positive real number such that
    $$\inf \R_{i\bar {i} j \bj}(g_0(t)) + \lambda>1.$$
Since $t_0<T$, we can choose $\lambda$ independent of $t$. 

We have
$$\p_{t}\Delta (\cdot) = \tr \left(g_0^{-1}(t)~g_0'(t)~g_0^{-1}(t)~\ddc
    (\cdot)\right).$$
Moreover, from (\ref{eta.etal}) it follows  that
$rg_0'(T)\le g_0(T)$, where $r$ is given by
    (\ref{def.r}) and $T=\log(r+1)$. Thus, there exists  a constant
    $C_0$, such that $g_0'(t)\le C_0~g_0(t)$ for any $t\in (0,T)$.
    In particular, since $g_0(t)+\ddc u>0$, we have
$$\begin{aligned}
(\p_t\Delta) u&= \tr
    \left(g_0^{-1}(t)~g_0'(t)~g_0^{-1}(t)~(g_0(t)+\ddc u)\right) - tr_{g_0(t)}g'_0(t)\\
&\le C_0n(n+\Delta u)-\tr_{g_0(t)}g'_0(t).
\end{aligned}$$ 
Thus, if $z=e^{-\lambda u} (n+\Delta u)$,
we get
\begin{equation}\label{eq.z}
\begin{aligned}
\p_t z - \Delta_{g(t)} z ~ \le ~ &
e^{-\lambda u} (n+\Delta f + n^2\inf \R_{i\bar{i} j \bj}(g_0(t)) +
      \Delta \log\det(g_0(t)/g_0) - \tr_{g_0(t)}g'_0(t))  \\
&+  (C_0 n - 1 + \lambda (n - \p_t u))z -  e^{\frac {\lambda u - F}{n+1}}(\inf \R_{i\bar {i} j \bj}(g_0(t))
    +\lambda)~
z^{\frac{n}{n-1}} \end{aligned}
\end{equation}

By the Maximum Principle and since, by prop. \ref{c0bound} and
    \ref{ut.bound}, $u$ and $u_t$ are bounded function in $(0,t_0)$,  the lemma follows.
\end{proof}

\vskip .5cm

From the above results, and from the Schauder estimate, we obtain a $C^1$-bound
    for the solution $u$. In fact, there exists a constant $C$ such that

$$\sup_{M\times [0,T)} \mid \nabla u\mid ~\le~ C~ (\sup_{M\times
    [0,T)}\mid \Delta u\mid +\sup_{M\times
    [0,T)}\mid  u\mid ) ~. $$

From the lemma and the fact that $g(t)=g_0(t)+\ddc u$ is a metric, we
    also obtain a $C^2$-bound for $u$.

\vskip .7 cm

Moreover, from prop. \ref{ut.bound} and from (\ref{scalar}), we get

\begin{lemma}\label{c2bound}
For any $t_0\in(0,T)$  there exist positive constants $C_0,C_1$, such that, for
    any $t<t_0$ and as long as there exists a solution for
    (\ref{scalar}),  we have:
$$C_0~g_0\le g(t)\le C_1~g_0.$$
\end{lemma}

\vskip .5 cm
Now, similarly to \cite{ca}, we can prove:

\begin{theorem}
Let $T=\log(r+1)$, with $r$ as in (\ref{def.r}).

There exists a bounded solution $u(x,t)\in
    C^\infty (M)$ for (\ref{scalar}), at  any time
    $t\in[0,T)$. Furthermore, since $u(x,t)$ are uniformly
    continuous in $[0,t)$, if $\lim _{t\to T} u(x,t)$ exists, it must be
    continuous.
\end{theorem}

\vskip 1 cm
\section{Singularitites of the Limiting Metric}

Let $M$ be a projective manifold with canonical class $K_M$ not nef,
    and let $g_0$ be a K\"ahler metric on $M$ belonging to the ample
    class $A$. 
In the previous section we have shown that
    as long as $A(t)$ is an
    ample class, i.e. as long as $t<T=\log(r+1)$ with $r$ as in
    (\ref{def.r}),
  there exists a solution $g(t)$ for the
    K\"ahler-Ricci flow (\ref{ricciflow}), that is in the class of
    $A(t)$.

On the other hand, at the finite time $T=\log(r+1)$, we have
$A(T)=\frac 1 {r+1}L$, which is a
    non-positive, but semi-ample, line bundle.
Thus, the associated pseudo-metric $g(T)$ will admit singularities,
    i.e. the set $$S=\{x\in M|\nexists ~ C_1,C_2>0 \text{ s.t. }~
    C_1g_0\le g(T) \le C_2 g_0 \text{ locally at }x\}$$ is not empty.

From now on, we will suppose that $M$ is of general type and
therefore
     $L$ is not only semi-ample, but also {\it big}, i.e., it defines a birational morphism
    $$c:M\to M'$$
    onto a projective variety $M'$. The exceptional
    locus $E$ of $c$, i.e. the subset of $M$ contracted by $c$, is the union
    of all the curves $C$, such that $L\cdot C=0$. Thus $E$ is the
    locus spanned by the extremal ray associated to $L$.

In general, we have

\begin{lemma}\label{contained}
The singular locus $S$ contains the exceptional set $E$. Moreover
    $g(T)$ is the pull-back of a (non-necessarily smooth) pseudo-metric on $M'$ of
    the form  $\alpha +dd^c u'$, where $\alpha$ is a smooth $(1,1)$-form and $u'$ is
    continuous (hence $L^1_{loc}$).
\end{lemma}

\begin{proof}
Since $g(T)$ belongs to the class of $\frac 1 {r+1} L$, it follows that for any
    curve $C$ in $M$,
$$\int_C g(T)=\frac {L\cdot C} {r+1}.$$

\noindent $g(T)$, being a limit of positive $(1,1)$-forms,
 it is  a
    non-negative $(1,1)-$form and therefore it follows that, for any curve $C$,
    such that $L\cdot C=0$, we have that $g(T)$ is zero along the
    directions
    tangential to $C$. In particular, $C$ is contained in $S$.

Since $L$ is semi-ample, there exists an ample class $A'$ on $M'$
such that $L=c^*A'$. Thus if we choose $\eta_L$ to be the pull
back of a non-negative $(1,1)$-form in $A'$, as in section
\ref{longtime}, we can write $g(t)=g_0(t)+\ddc u$, with
$$g_0(t)=\frac 1 r (a(t)\eta_L+ b(t)g_0).$$
Since both $g_0(T)=\frac 1 {r+1} \eta_L$ and $g(T)$ are zero along
the direction tangential to $C$, it follows that $\ddc u=0$ along
$C$ and therefore $ \Delta _C u=0$ (here $\Delta _C$
is the complex Laplacian on $C$). Thus, since $C$ is compact, $u$
is smooth along $C$ (by {\it elliptic regularity}) and it induces
a function $u'$ on $M'$. Since $u$ is continuous, it must be the
case that also $u'$ is. Thus $g(T)=\eta_L+c^*(\ddc u')$ is the
pull-back of a $(1,1)$ form on $M'$ with the required properties.
\end{proof}

\begin{remark}
Note that this lemma holds also in the case that $L$ is only
semi-ample (i.e., without requiring that $K_M$ be big), in that
case though, one would get a map $c: M\to S$ with $\dim(S)<\dim (M)$
and $g(T)$ would be degenerate along the fibers.
\end{remark}
\vskip .5 cm

We will need the following:

\begin{lemma}\label{bis.estimate}
Let $f:M \to N$ be a smooth holomorphic map of K\"ahler manifolds,
and let $g=g_N$ be a K\"ahler metric on $N$ such that
$h=f^*g_N\geq0$. Then $\R_{i\bar i j\bar j}(h)$ is bounded.
\end{lemma}

\begin{proof}
The boundedness of the bisectional curvature being a local matter,
we can pick a point $p\in M$ and show the boundedness around that
point.

Since the property of being bounded is independent of the
coordinates chosen, we can choose coordinates $\{z_1, \cdots ,
z_n\}$ centered around $p$ and coordinates $\{ w_1, \cdots ,
w_n\}$ centered around $f(p)$ such that $h=\sum _i \lambda _i
dz_i\wedge dz_{\bar i}$ and $g= \sum _i \mu _i dw_i \wedge
dw_{\bar i}$. In particular:

$$\lambda _k = \sum _{i=1}^n \mu _i \frac{\p f_i}{\p z_k} \overline {\frac{\p f_i}{\p
z_k}}$$ \noindent and since $g$ is positive, this vanishes exactly
if and only if $ | \frac{\p f_i}{\p z_k}  |^2=\frac{\p f _i}{\p
z_k} \overline{\frac{\p f_i}{\p z_k}}=0$ for every $i$. Let $i_o$
be such that $| \frac{\p f_{i_o}}{\p z_k}|^2$ has minimal order of
vanishing at $p$ (this is the order of vanishing of $\lambda _k$)

One computes:
$$ \frac{\p \lambda _k}{\p z_j} =\sum _{i=1}^n \frac{\p \mu _k}{\p
z_j} \frac{\p f_i}{\p z_k} \overline {\frac{\p f_i}{\p z_k}} +
\sum _{i=1}^n \mu _k \frac{\p ^2 f_i}{\p z_j\p z_k} \overline
{\frac{\p f_i}{\p z_k}}$$ and the analogous formula for $\frac{\p
\lambda _k}{\p \bar z_j}$:

$$  \frac{\p \lambda _k}{\p \bar z_j} =\sum _{i=1}^n \frac{\p \mu _k}{\p
\bar z_j} \frac{\p f_i}{\p z_k} \overline {\frac{\p f_i}{\p z_k}}
+ \sum _{i=1}^n \mu _k \frac{\p  f_i}{\p z_k} \overline {\frac{\p
^2 f_i}{\p z_h\p z_k}} .$$

\noindent
Therefore, $\frac{\p \lambda _k}{\p z_j}=| \frac{\p
f_{i_o}}{\p z_k}|^2 A_1 + \frac{\p}{\p z_j}\left( | \frac{\p
f_{i_o}}{\p z_k}|^2 \right) B_1$ and $\frac{\p \lambda _k}{\p \bar
z_j}= | \frac{\p f_{i_o}}{\p z_k}|^2 A_2 + \frac{\p}{\p \bar
z_j}\left( | \frac{\p f_{i_o}}{\p z_k}|^2 \right)B_2$

\noindent From this it is easy to see that $\frac{\p \lambda
_k}{\p z_j}\frac{\p \lambda _k}{\p \bar z_j}$ is divisible by $|
\frac{\p f_{i_o}}{\p z_k}|^2$.

Since in these coordinates the bisectional curvature equals:
$$\R_{k\bar k j\bar j}(h)= -\frac{\p ^2 \lambda _k}{\p z_j \p\bar
z_j}+ \frac{1}{\lambda _k} \sum _ k \frac{\p \lambda _k}{\p z_j}
\frac{\p \lambda _k}{\p \bar z_j},$$ we confirm the conclusion of
the theorem.
\end{proof}
\vskip .5 cm

In order to have a better understanding of the singularities of
the pseudo-metric induced on $M'$, we will first study the locus:
$$S_0=\{x\in M|\det(g(T))=0\}.$$
Obviously $S_0$ is contained in the singular locus $S$ for $g(T)$.

Although not quite standard, we will call a birational contraction
$c:M\to
    M'$ {\em divisorial} if the exceptional locus $E$ of $c$ is of
    pure codimension $1$ and
    $M'$ is $\QQ-$factorial, i.e. for any Weil divisor $D$ on $M'$,
    there exists an integer $m$ such that $mD$ can be locally defined
    by one equation. This last requirement is superfluous if $c$ is an
    extremal contraction. Moreover in the case of a projective
    surface $S$, the contraction of any $K_S$-negative curve is
    divisorial.

\begin{proposition}\label{locus.s0}
The  locus $S_0$ where $\det g(T)=0$ is contained in a
    proper sub-variety of $M$. Moreover, if $c:M\to M'$ is a divisorial
    contraction, then the  locus $S_0$  coincides with
    the exceptional locus $E$ (when viewed as sets).
\end{proposition}

\begin{proof}

We will first show that $S_0\subset F$ for some analytic
subvariety $F$. Let $A=A(0)$ be the class of the initial metric
$g_0$. Since $K_M$ is big, by Riemann-Roch it follows that there
    exists a positive integer $m$ and a $\QQ-$effective $\QQ-$divisor
    $D$ (i.e. $dD$ is an effective integral divisor for some positive
    integer $d$)
   such that $mK_M=A+D.$

By lemma \ref{rescale}, we can rescale $g_0$ by $1/m$, so that we can
    suppose without loss of generality, that $m=1$ and $K_M=A+D$.

Let us choose $\eta_0=-g_0-\Ric_{g_0}$. For this choice, in the class
    of $K_M-A$,  the equation
(\ref{scalar}) becomes:
\begin{equation}\label{eta0}
\p_{t} u= \log ~ \det\left( \frac{g_0+a(t)\eta_0+\ddc u}{g_0} \right)-u
\end{equation}

Therefore, by differentiating with respect to $t$, we have, as in
(\ref{heat.scalar}), that
    $v=\partial_t u$ is a solution of:
\begin{equation}\label{heat}
\left\{
\begin{aligned}
&\partial_t v = \Delta_{g(t)} v + a'(t) \tr_{g(t)} \eta_0 - v\\
&v(x,0)=0
\end{aligned}\right.
\end{equation}

From (\ref{eta0}) and since by proposition \ref{c0bound} the
function $u$ is
    uniformly bounded in $(0,T)$, it follows immediately that the locus $S_0$ where
$\det g(T)=0$, coincides with the locus where $\displaystyle \lim_{t\to T}v= -\infty$.

Since $\eta_0$ is in the class of $c_1(D)$, and $D$ is
$\QQ-$effective,
    there exists a positive integer $d$, such
    that $d\eta_0$ is in the class of the effective divisor $dD$.

Let $\sigma$ be a non zero holomorphic section in $\HH^0(M,dD)$ with
    zero locus $F$, and
    let $h$ be any hermitian metric on $dD$. Then, by Chern theory, there exists a
    smooth function $\phi_0$ on $M$ such that $\frac {\sqrt {-1}}{2\pi} \Theta_h(D)=\eta_0 +
    \ddc \phi_0$.

Moreover, by the Lelong-Poincar\'e formula (e.g.
    see  (3.11) in \cite{D})
we have
$$\ddc \log (\|\sigma\|^2_h)=[F] -  \frac {\sqrt {-1}}{2\pi} \Theta_h(dD)$$
where the equality is understood in the sense of currents and $[F]$ denotes the current of integration associated to
    $F$, defined by, for any $(n-1,n-1)-$form $\alpha$,
$$<[F],\alpha>=\int_F\alpha.$$

Therefore if $\phi=\phi_0+\frac 1 d\log(\|\sigma\|^2_h)$ on $M\setminus
    F$, we have
\begin{equation}\label{lp}
\eta_0=\frac 1 d [F]-\ddc \phi
\end{equation}

\vskip .5 cm

Let $\overline v= v-a'(t)\phi$, then $\overline v(\cdot, t)$ is a smooth function
    on $M\setminus F$, defined for any $t<T$. Moreover, since
    $a''(t)=-a'(t)$, we have that
$\overline v$ satisfies the equation:

$$
\begin{aligned}
\p_t \overline v & = \p_t v - a''(t)\phi \\
        & =\Delta_{g(t)} v + a'(t) \tr_{g(t)} \eta_0 - v + a'(t)\phi\\
        & = \Delta_{g(t)} \overline v + a'(t) \tr_{g(t)}\ddc\phi +
    a'(t)\tr_{g(t)} \eta_0 - \overline v
\end{aligned}
$$

Thus by (\ref{lp}), it follows  that  $\overline v$ is a
    solution of
\begin{equation}
\label{heat.mod}
\left \{
\begin{aligned}
\overline v_t & = \Delta_{g(t)} \overline v - \overline v \qquad &\text{on }
    (M\setminus F)\times (0,T)\\
\overline v&=-\phi \qquad & \text{on }(M\setminus F)\times \{t=0\}
\end{aligned}
\right .
\end{equation}

Since for any $t\in [0,T)$ and for any $x$ approaching $F$,
we have that $\overline v(x,t)\to +\infty$, then the function $\overline v$ must
    admit a minimum inside $M\setminus F$ for any time $t$.
Therefore by the Minimum Principle, applied at (\ref{heat.mod}), it follows that
$$\overline v(t)\ge \inf\overline v(\cdot, 0)= - \max \phi >-\infty.$$

Therefore $v = \overline v + a'(t)\phi \ge -\max \phi + a'(t)\phi$, and in
    particular the locus $S_0$ where $v\to -\infty$ must be contained inside
    $F$ (since $a'(t)= e_t$ is bounded in $t$, for $t\geq 0$).

\vskip .5 cm

Let us now prove that under the assumption that $c:M\to M'$ is a
divisorial contraction, then the locus $S_o$ coincides with the
    exceptional locus $E$ of $c$. Let $E_i$ be the irreducible components
    of $E$. Then, since we are assuming that $M'$ is $\QQ-$factorial,  the relative Picard group of $c$ is generated by
    $E_i$, i.e. $$\Pic(M)/c^*\Pic(M')=<E_i>.$$
Therefore $A=c^* B - \sum \delta_i E_i$, for some ample class $B$
    on $M'$ and positive constants $\delta_i$.

If $\alpha_i$ is the discrepancy of $E_i$ with respect to the morphism
    $c$, i.e. if
$$\alpha_i = \mult_{E_i}(K_{M/M'}),$$
then, since $M'$ has terminal singularities (e.g. see \cite{kmm}), it
    follows that $\alpha_i>0$.

Thus
$$L=rK_M + A = c^*(rK_{M'}+B) + \sum (r\alpha_i-\delta_i)E_i, $$
from which it follows that $r\alpha_i-\delta_i=0$.

\noindent
 Therefore one has:
    $$K_M-A=c^*(K_{M'}-B)+(r+1)\sum \alpha_i E_i.$$
\noindent
From lemma \ref{contained}, it follows that $g(T)$ is the pull-back of
a non-negative $(1,1)$-form $g'$  in an ample class of $M'$,
    therefore, it follows from lemma \ref{trace} that we can choose a $(1,1)$-form
    $\eta'$ in the class of $K_{M'}-B$, such that
\begin{equation}\label{tracebound}
\tr_{g'}\eta'>C,
\end{equation}
    for some constant $C$.

Thus, similarly to (\ref{lp}), there exists a smooth function
$\phi_0$,
    such that if $$\phi=\phi_0+(r+1)\sum \alpha_i \log(\|\sigma_i\|^2_{h_i}),$$
where
    $\sigma_i\in \HH^0(M,E_i)$ are non-zero holomorphic sections and
    $h_i$ are
    hermitian metrics on $E_i$, then
$$\eta_0=(r+1)\sum \alpha_i[E_i] + c^*\eta' - \ddc \phi.$$
\noindent
Thus, from (\ref{tracebound}), equation (\ref{heat.mod}) becomes
$$
\left \{
\begin{aligned}
\overline v_t & \ge \Delta_{g(t)} \overline v - \overline v +C \qquad &\text{on }
    (M\setminus E)\times (0,T)\\
\overline v&=-\phi \qquad & \text{on }(M\setminus E)\times \{t=0\}
\end{aligned}
\right .
$$
and similarly to what we did earlier, it follows that there exists
a lower bound
    for $\overline v$. Thus, since $v=\overline v + a'(t)\phi$, and
    $T=\log(r+1)$, there exists a constant $B$ such that
$$\begin{aligned}
{\p_t u}_{|_{t=T}}=v_{|_{t=T}}&\ge B + a'(T)(r+1)\sum \alpha_i \log(\|\sigma_i\|^2_{h_i})=\\
&=B+\sum \alpha_i \log(\|\sigma_i\|^{2}_{h_i})
\end{aligned}
$$
from which it follows that $v$ is bounded outside $E$.

Therefore, since by proposition \ref{c0bound} the
    function $u$ is bounded, eq. (\ref{eta0}) yields the existence of a positive constant $B'$ such that
\begin{equation}\label{det.comparison}
\det(g(T)/g_0)\ge B' \prod\|\sigma_i\|^{2\alpha_i}_{h_i}.
\end{equation}
In particular, the locus $S_0$ where $\det g(T)=0$ coincides with
    $E$.
\end{proof}

\vskip .5 cm
We are now ready to prove the main result of this section
\begin{theorem}\label{smooth}
If $c:M\to M'$ is a divisorial contraction and $M'$ is smooth, then
    the singular locus $S$ for $g(T)$ coincides with the exceptional locus
    $E$. Moreover $g(T)$ induces a smooth metric on $M'$.
\end{theorem}

\vskip .5 cm

We will need the following:
\begin{lemma}\label{laplacian.comparison}
    If the bisectional curvature $\R_{i \bar i j \bar j}$ of $g_0(t)$ is uniformly bounded in
    $[0,T)$, then there exists a constant $C$ such that
    $$\Delta_{g_0(t)}u + \log(\det(g_0(t)/g_0))<C.$$
\end{lemma}
\begin{proof}
Let
$$\overline z=e^{-\lambda
    u}(n+\Delta_{g_0(t)}u+\log \det(g_0(t)/g_0)+f).$$
where, as in the proof of lemma \ref{laplacianestimate}, $\lambda$ is
    a postitive constant  such that 
$$\inf \R_{i\bar i j \bar j}(g_0(t))+\lambda >1.$$

Then, by equation (\ref{eq.z}), there exist constant $C_0,C_1,C_2$ with $C_2<0$ and such that
$$\p_t \overline z - \Delta_{g(t)} \overline z ~ \le
e^{-\lambda u} (n + n^2\inf \R_{i\bar i j \bar j}(g_0(t))+C_0) +C_1
    \overline z + C_2 \overline z^{\frac n {n-1}}  $$
Thus, the claim follows immediately by the Maximum Principle.
\end{proof}

\vskip .5 cm

\begin{proof}[Proof of Theorem \ref{smooth}]
\

Let $\eta_L$ be induced by a smooth metric on $M'$. From
    lemma \ref{bis.estimate}, it follows that
    the bisectional curvature of $g_0(t)$ is uniformly bounded.
Thus by lemma \ref{laplacian.comparison}, the metric $g(t)$
    remains bounded outside the exceptional locus of $c:M\to M'$.

Moreover,
    there exists a smooth function $\psi$ on $M$,
    such that,
using the same notation as in the proof of proposition \ref{locus.s0} , we have
\begin{equation}\label{det.etal}
\det(\eta_L/g_0)= e^\psi \prod\|\sigma_i\|^{2\alpha_i}_{h_i}.
\end{equation}
\noindent
In particular by (\ref{det.comparison}), it follows that there exists
    a positive constant $B_0$ such that
$$\det(g(t)/g_0(t))>B_0,$$
for any $t\in (0,T)$.
\noindent
We now want to show that there exists an upper bound for
    $\det(g(t)/g_0(t))$.
\noindent
Let $\phi =e^t\log \det (g_0(t)/g_0)$ and $w=e^t \p_t u - \phi$.
Then, equation (\ref{heat.scalar}) becomes
\begin{equation}\label{eq.w}
\left\{
\begin{aligned}
&\partial_t w = \Delta_{g(t)} w + \tr_{g(t)} (\eta + \ddc \phi) - \phi'\\
&w(x,0)=f
\end{aligned}\right.
\end{equation}
where $\eta$ depends on $\eta_L$, as in  (\ref{eta.etal}).
\noindent
Let us study the term $\Phi= \tr_{g(t)} (\eta + \ddc \phi) - \phi'$ in
    equation (\ref{eq.w}).
We have:
$$
\Phi = \tr [(g^{-1}(t)-g_0^{-1}(t))\eta] + \tr_{g(t)}\ddc\phi - \phi.$$
 From (\ref{det.etal}) we have
$$\ddc \phi_{|_{t=T}} = (r+1)\sum \alpha_i\ddc \log ||\sigma_i||_{h_i}^2$$
in fact, we have that there exists a $\delta>0$ and constants $A$ and
$B$ such that for every
$t\in (T-\delta, T)$, one has:
$$\ddc \phi  \leq Ae^t \sum \alpha_i\ddc \log ||\sigma_i||_{h_i}^2$$
and
$$\phi  \leq Be^t \sum \alpha_i\ddc ||\sigma_i||_{h_i}^2$$
from which it easely follows
that every term of $\psi$ is integrable with respect to $t$ in
    $(0,T)$, and by the Maximum Principle it follows that the solution $w$ for
    (\ref{eq.w}) is bounded. In particular $\det(g(t)/g_0(t))$ is
    uniformly bounded.
\end{proof}
\vskip 1 cm
\section {Maximal existence time. Case II:  $K_M$ nef}\label{longtime2}

In this section, we consider the case of a smooth projective variety
    $M$ such that the canonical class $K_M$ is  nef. As above, let
    $A$ be the ample class that represents the
    initial metric $g_0$.

By prop. \ref{class.eq} the solution $g(t)$ for
    the K\"ahler-Ricci flow (\ref{ricciflow}) belongs to the class
$$A(t)=A+a(t)(K_M-A),$$
where $a(t)=1-e^{-t}$.
Since $K_M$ is nef, it follows that $A(t)$ is ample for any $t>0$.

Moreover, by the base point free theorem, we have that $L=K_M$ is
    semi-ample and therefore, as in \ref{eta.etal}, we can write
    \begin{equation}\label{g0.nef}
    g_0=a(t)\eta_L+b(t)g_0,
    \end{equation}
where $\eta_L$ is a non-negative
    $(1,1)$-form in the class of $K_M$ and $b(t)=e^{-t}$.

Thus the solution $g(t)$ can be written as $g(t)=g_0(t)+\ddc u,$ for some
    function $u$ that satisfies equation (\ref{scalar}).

\begin{proposition}\label{case.nef}
If $M$ is a projective variety with nef canonical line bundle, then
    the K\"ahler-Ricci flow (\ref{ricciflow}) admits a smooth solution $g(t)$
    for any time $t\in(0,+\infty)$.

 \end{proposition}

\begin{proof}
As in lemma \ref{supersolution}, we are going to show the existence of
    a bounded super-solution $u^+$ and a bounded sub-solution $u^-$
    for (\ref{scalar}).

Since $a(t)$ and $b(t)$ are bounded and $M$ is compact, there exists a
    constant $K$, such that $\log\det(g_0(t))+f< K$. Therefore, as in
    lemma \ref{supersolution}, in order to find a super-solution
    $u^+$, it is enough to solve $\partial_t {u^+} + u^+=K$, with $u^+(0)=0$.

Similarly, for any  $t_0<+\infty$, it is easy to find a bounded and space-independent
    sub-solution $u^-$ (depending on $t_0$), defined for any $t\in
    (0,t_0)$. 

Moreover, as in proposition \ref{ut.bound}, $v=\p_{t} u$  is a
    solution for
\begin{equation}\label{heat.scalar.nef}
\left\{
\begin{aligned}
\partial_t v &= \Delta_{g(t)} v + a'(t) \tr_{g(t)} (\eta_L-g_0 ) - v\\
v(0) &= f
\end{aligned}\right.
\end{equation}

By lemma \ref{trace}, for any $t_0<+\infty$, we can modify the equation
    (\ref{heat.scalar.nef}), by choosing a $(1,1)-$form $\eta'$
    (resp. $\eta''$) cohomologous to $\eta_L-g_0$ and such that
    $\tr_{g(t)}\eta' > C$ (resp. $\tr_{g(t)}< C$) for some constant
    $C$. From that, it follows that  $v$ is bounded in the interval $(0,t_0)$.

Following exactly the same lines as in section \ref{longtime}, we
    obtain the same $C^1$ and $C^2$ estimates for $u$. Thus, for any
    $t_0\in(0,+\infty)$,  there
    exists a bounded solution $u$ for (\ref{scalar}) in $(0,t_0)$. 

Thus, the claim follows. 
\end{proof}

\vskip .5 cm

From the  proof of the previous proposition, it follows  that the solution $u$
    for (\ref{scalar}) is uniformly bounded from above for any $t\in (0,+\infty)$.

As in the previous section, we will suppose from now on, that $K_M$ is also
    big. Since by the nefness assumption, $K_M$ is  semi-ample, it
    defines a birational morphism $c:M\to X$, onto the canonical model $X$ for
    $M$. By a classic result in algebraic geometry, $X$ has canonical
    singularities, and $K_X$ is ample (e.g. see \cite{kmm}).

In order to study the singularities for the limit metric for the
    K\"ahler-Ricci flow in this situation, we are going to restrict
    our self to the case $X$ admits only orbifold singularities. 
 In particular,
    this always holds if $M$ is a surface of general type.

Thus, we have

\begin{proposition}
Let $M$ be a projective manifold of general type with nef canonical
    bundle and such that its canonical model admits only orbifold
    singularitites.

Then equation (\ref{scalar})
admits a uniformly bounded solution for any $t\in (0,+\infty)$.
\end{proposition}

\begin{proof}
By the proof of proposition \ref{case.nef}, and by the comparison
    principle, it is enough to show that
    there exists a bounded sub-solution in some interval $(t_0,+\infty)$.

Since $K_M$ is semi-ample, there exists a smooth non-negative $(1,1)$-form
    $\eta_L$ in the class of $K_M$. Moreover, since by assumption $X$
    is an orbifold, we can choose $\tilde \eta_L$ as the pull-back of
    an orbifold metric on $X$ (e.g. see \cite{DK}).
In particular, if $c:M\to X$ is the canonical map defined by $K_M$,
    since $K_M=c^*K_X$,  there exists a smooth function $\psi$, such
    that
\begin{equation}\label{orb.volume}
\det  {\tilde \eta_L} =    \det g_0\cdot e^{-\psi}.
\end{equation}

Moreover, there exists a continuous function $\phi_0$ (that is smooth in
    the orbifold sense), such that
$$\tilde \eta_L=\eta_L + \ddc \phi_0.$$
It is possible to approximate $\phi_0$ by a sequence of smooth (in the
    ordinary sense) function on $M$ (e.g. see \cite{B}). Therefore we
    can choose $\phi(t)\in C^\infty(M)$ uniformly bounded and such
    that
$$\lim_{t\to\infty} \phi(t)=\phi_0 \qquad \text{and} \qquad \lim_{t\to\infty} \ddc\phi(t)=\ddc\phi_0  .$$

For any $\rho\in (0,1)$, let us define
$$G_{\rho}:=\log\det\rho\cdot\left(\frac{g_0(t)+\ddc \phi(t)}{g_0} \right) + \psi$$
with $g_0(t)$ as in   (\ref{g0.nef}). Then, by (\ref{orb.volume}), it follows
    that  $\lim_{t\to \infty} G_\rho=n\log \rho$ is a constant.

Fixed a constant $$ C\ge\sup_{M\times(0,+\infty)}
    (\phi(t)+\phi'(t)-f),$$ let
$$u^-:=e^{-t}\cdot\int_{0}^t e^s G_\rho~ds + \phi(t) -C.$$

Clearly $u^-$ is bounded.
We want to show that if $t_0$ is sufficiently large (depending on
    $\rho$), then
    $u^-$ is a sub-solution for (\ref{scalar}) in the
    interval $(t_0,+\infty)$. This will follow from:

\vskip .5 cm
{\bf Claim:} There exists $t_0$ such that
\begin{equation}\label{claim}
\det(\rho (g_0(t) + \ddc \phi(t)))\le \det (g_0(t) + \ddc u^-)
\qquad \text{for }t\in(t_0,+\infty)
\end{equation}

\vskip .5 cm

In fact, we have

$$\begin{aligned}
  \p_t u^- + u^- = &~ G_\rho +\phi'(t) + \phi(t) - C  \\
    \le ~& \log\det \rho \left(\frac{ g_0(t)+\ddc \phi(t)}{g_0}\right) + f\\
   \le ~& \log\det\left(\frac{g_0(t)+ \ddc u^-}{g_0}\right) + f.
\end{aligned}
$$

where the second line follows from the choice of $C$, while the third
    line follows from (\ref{claim}).
Thus $u^-$ is a sub-solution for (\ref{scalar}), and in particular
there exists a constant $C'$ such that    $u\ge u^-+C'$, for any $t\in (t_0,+\infty)$.

\vskip .5 cm

In order to prove the claim, it is enough to observe that by taking
    the limit $t\to +\infty$, we have
$$\begin{aligned}
\lim_{t\to +\infty}\det(\rho (g_0(t) + \ddc \phi(t)))
&=\det(\rho (\eta_L + \ddc \phi_0)) \\
& <\det(\tilde \eta_L)\\
& = \lim_{t\to +\infty} \det (g_0(t) + \ddc u^-).
\end{aligned}
$$
Thus, by the compactness of $M$,  the claim follows.
\end{proof}

\end{document}